\documentclass[12pt]{article}
   \def\sqr#1#2{$\vcenter{\hrule height.#2pt
   \hbox{\vrule width.#2pt height#1pt \kern#1pt
   \vrule width.#2pt}
   \hrule height.#2pt}$}
   \def\qed{\sqr53}
   \newfont{\eightrm}{cmr8}
   \newfont{\ninerm}{cmr9}
   \newfont{\nineit}{cmti9}
   \newfont{\eightit}{cmti8}

\begin{document}
\begin{center}
\vskip 0.6cm {\large\bf Fundamental Cycles and Graph Embeddings
\footnote{Supported by Natural Science Foundation of China ( Under
the Granted Number 10271048, 10671073)}}

\vskip 0.2cm

Ren Han\footnote {This work is partially supported by Science and
Technology Commission of Shanghai Municipality (07XD14011) and
Shanghai Leading Academic Discipline Project, Project
Number£ºB407},\,Zhao Hongtao and Li Haoling

\vskip 0.3cm

 {\footnotesize{\it Dept. of Mathematics, East China Normal University}}\\[5pt]
 {\footnotesize{\it Shanghai 200062, P.R.China}}\\[5pt]
{\footnotesize{\it E-mail: hren@math.ecnu.edu.cn}}

\vskip .3cm
\end{center}

 \vskip 0.3cm
 \noindent {\bf Abstract:}
 \normalsize
 In this paper, we investigate fundamental cycles in a graph~$G$~and their
 relations with graph embeddings.~We show
 that a graph~$G $~may be embedded in an
 orientable surface with genus at least $g$ if and only if for any spanning tree $T$,
 there exists a sequence of
 fundamental cycles $C_1,~C_2,~\cdots,C_{2g}$ with $C_{2i-1}\cap C_{2i}\neq\phi$~
 for~$1\leq~i\leq g $. In particular, among $\beta(G)$
 fundamental cycles of any spanning tree $T$ of a graph $G$, there
 are exactly $2\gamma_M(G)$ cycles
 $C_1,~C_2,~\cdots,C_{2\gamma_M(G)}$ such that $C_{2i-1}\cap C_{2i}\neq\phi$~
 for~$1\leq~i\leq\gamma_M(G)$, where $\beta(G)$ and $\gamma_M(G)$ are, respectively, the Betti
 number and the maximum genus of $G$. This implies that it is possible to construct an
 orientable embedding with large genus of a graph $G$ from an arbitrary spanning tree $T$( which
 may have very large number of odd components in $G\backslash E(T)$).
 This is different from the earlier work of Xuong and Liu[9,6],
 where spanning trees with small odd components are needed. In fact,
 this makes a common generalization of Xuong[9],Liu[6] and Fu et al[2].
 Further more, we show that (1).This result is useful in locating the
 maximum genus of a graph having a specific edge-cut.
 Some known results for embedded graphs are also concluded;(2).The maximum
 genus problem may be reduced to the maximum matching problem.
 Based on this result and the algorithm of
 Micali-Vazirani[8], we present a new efficient algorithm to determine the
 maximum genus of a graph in $O({(\beta(G))}^{\frac{5}{2}})$ steps.
 Our method is straight and quite deferent from the algorithm of Furst,Gross
 and McGeoch[3] which depends on a result of Giles[4]where matroid parity method is needed.\\
\noindent{\bf Keyword~:}   Fundamental cycles, Maximum
genus,upper-embedded~.\\
\noindent {\bf AMS 2000:}\, Primary 05C10, secondary 05C70

\newpage

\noindent\section{ Definitions and Notations}

 The graph considered here is finite and undirected and, furthermore,
is connected unless it is stated otherwise. In general, multiple
edges and loops are allowed. Terminology and notation without
explicit explanation follows as from [1,6,7].

\vskip 0.2cm

 By a {\it surface}, denoted by $S$, we mean a compact and connected
2-manifold without boundary. It is well known from elementary
topology that surfaces can be divided into two classes:~{\it
orientable}~and {\it nonorientable} ones. An {\it orientable
surface} can be viewed as a sphere attached $h$ handles, while a
{\it nonorientable surface} as a sphere attached $k$ crosscaps. The
number $h$ or $k$ is called the {\it genus} of the surface. A {\it
cellular embedding} of a graph $G$ into a surface $S$ is a
continuous one-to-one mapping $\phi$:~$G\rightarrow S$ such that
each component of $G\backslash\phi(G)$ is homeomorphic to an open
disc, called a {\it face} of $G$~(~with respect to this embedding
$\phi$) and $\phi$ is called a {\it cellular embedding}( or {\it
embedding }as some scholars called ). A cycle ( curve ) $C$ in an
embedded graph in a surface $\sum$ is called {\it surface separating
} if $\sum-C$ is disconnected. In particular, if $\sum-C$ has an
open disc, denoted by int$(C)$, then $C$ is called {\it
contractible} (otherwise, $C$ is {\it noncontractible}), and
int$(C)+C=$Int$(C)$ is the {\it inner part} of $C$ . The other part
of $\sum-C$ is called {\it exterior of} $C$ and is denoted by
Ext$(C)$.

\vskip 0.2cm
 Recall that the {\it maximum genus} $\gamma_{M}(G)$ of a graph $G$ is the largest
integer $k$ such that $G$ has an embedding in an orientable surface
with genus $k$. Since any graph $G$ embedded in a surface has at
least one face, Euler's formula shows that
$\gamma_{M}(G)\leq\lfloor\frac{\beta(G)}{2}\rfloor$, where
$\beta(G)=|E(G)|-|V(G)|+1$ is known as {\it Betti number} of $G$
(~which is equal to the cyclic number of $G$). A graph is {\it
upper-embeddable}~if~$\gamma_{M}(G)=\lfloor\frac{\beta(G)}{2}\rfloor$.

\vskip 0.2cm

 Let $G$ be a graph and $T$ be a spanning tree of $G$.
It is clear that for any edge $e\in {E(G)\backslash E(T)}$, $T+e$
contains a unique cycle of $G$, denoted by $C_T(e)$, which is called
a {\it fundamental cycle} of $G$ (~with respect to the spanning tree
$T$ of $G$). If a pair of edges $e_1$ and $e_2$ have a common end
vertex in a graph $G$, then we say that the pair~$\langle e_1, e_2
\rangle$~ is an {\it adjacent-edge pair} in $G$.Let $G_1$ and $G_2$
be a pair of disjoint subgraphs of $G$. Then $E [G_1 , G_2 ]$ is the
set of edges with their ends in $G_1$ and $G_2$,respectively.

\vskip 0.2cm

Denote by $\xi(G,T)$ the number of components of $G\backslash E(T)$
with an odd number of edges. Then the Betti deficiency of $G$
denoted by $\xi(G)$ is defined as the value $min_T\xi(G,T)$, where
the minimum is taken over all spanning trees $T$ of $G$. A spanning
tree $T$ of $G$ is said to be an {\it optimal spanning tree} if
$\xi(G,T)=\xi(G)$.

\vskip 0.2cm

 In the following, the
paper is organized as follows: in {\S 2} we give a good
characterization (i.e., Theorems 1 and 2) of maximum genus; {\S 3}
will concentrate on the applications of Theorems 1 and 2 and their
refined form; {\S 4} will show that finding the maximum genus of a
graph $G$ is, in some extend, equivalent to the problem of finding a
maximum matching in a specific graph $G_M$ called {\it the
fundamental intersecting graph} of $G$ and presents an efficient
algorithm in finding the maximum genus of a graph.
\section{A Good Characterization }
%\newtheorem{lemma} {Lemma}
%\begin{lemma}{\upshape{[6,9]}}
\noindent{\bf Lemma 1[6,9]} {\it Let $G$ be a graph, then\\
(1)\,$\gamma_{M}(G)=\frac{1}{2}(\beta(G)-\xi(G))$;\\
(2)\, $G$ is upper embeddable if and only if $\xi(G)\leq1$.}
%\end{lemma}
%\newtheorem{theorem}{Theorem}
%\begin{theorem}
\vskip 0.2cm

 \noindent{\bf Theorem 1.}{\it If a graph $G$ contains a
spanning tree $T$ such that there exist $2g$ fundamental cycles
$C_1, C_2, \cdots, C_{2g}$ with $C_{2i-1 }\cap C_{2i-1}\neq\phi$,
for $i=1, 2\cdots g$, then $G$ may be embedded in an orientable
surface with genus at least $g$.}

\vskip 0.2cm

\noindent{\bf Proof} Let~$G$~and~$T$~be as assumed and~$e_1, e_2,
\cdots, e_{2g}$~be edges in $E(G)\setminus E(T)$~such that $C_i$~is
the unique cycle in $T+e_i(1\leq i\leq 2g).$ We may suppose further
that $G=T+\{e_1, e_2 , \cdots, e_{2g}\}$ by Xuong's constructive
proof of maximum genus formula[9].~Let~$G_0=T$,~and~$G_1=G_0+\{e_1,
e_2\}$. Then we have the following.

\vskip 0.2cm

\noindent{\bf Claim 1.}
 ~~$\xi(G_1)\leq\xi(G_0)$~$(
\Leftrightarrow \gamma_{M}(G_1)\geq\gamma_{M}(G_0)+1).$
 \vskip 0.2cm
 To see this, we observe
that~$\beta(G_0)=\beta(G_1)-2$,~and
so,~$\xi(G_1)\equiv\xi(G_0)(~mod~2~)$.~\\
If~$\xi(G_1)\geq\xi(G_0)~+~2$, then we have one of the following
situations:\\
\noindent
 (1).~Both $e_1$~and $e_2$~have, respectively, their ends in distinct even components in $E(G_0)\backslash E(T)$
 (~As shown in left hand side of Fig.1).\\
\noindent
(2).~Both $e_1$~and $e_2$~have, respectively, their ends
in the same even components in $E(G_0)\backslash E(T)$
 (~As shown in center of Fig.1).\\
\noindent
(3).~Exactly, one of $e_1$~and $e_2$, say $e_1$, joins two
even components of $E(G_0)\backslash
 E(T)$, while~$e_2$~has two ends in the same even components in $E(G_0)\backslash
 E(T)$~(~As shown in right hand side of Fig.1).

%\vskip 0.2cm

 \begin{figure}[h]
\begin{picture}(220,100)(0,0)
\thicklines \put(30,20){\circle{40}} \put(40,20){\circle{1.5}}
\put(90,20){\circle{40}}\put(80,20){\circle{1.5}}
\put(30,80){\circle{40}} \put(40,80){\circle{1.5}}
\put(90,80){\circle{40}}\put(80,80){\circle{1.5}}
\put(150,70){\circle{40}}\put(210,70){\circle{40}}
\put(140,60){\circle{1.5}}\put(160,60){\circle{1.5}}
\put(200,60){\circle{1.5}}\put(220,60){\circle{1.5}}
\put(320,70){\circle{1.5}} \put(280,70){\circle{1.5}}
\put(330,70){\circle{40}}\put(270,70){\circle{40}}
\put(390,70){\circle{40}}
\put(380,60){\circle{1.5}}\put(400,60){\circle{1.5}}
 \qbezier(140,60)(150,15)(160,60)\qbezier(380,60)(390,15)(400,60)
\qbezier(200,60)(210,15)(220,60) \thinlines
\put(25,30){\line(2,1){10}}\put(75,30){\line(2,1){10}}
\put(40,20){\line(2,0){40}}\put(25,70){\line(2,1){10}}
\put(40,80){\line(2,0){40}}\put(75,90){\line(2,1){10}}
 \put(15,75){\line(2,1){25}}\put(76,80){\line(2,1){25}}
\put(25,70){\line(2,1){20}}\put(80,70){\line(2,1){20}}
\put(35,65){\line(2,1){10}}\put(90,65){\line(2,1){10}}
 \put(15,15){\line(2,1){25}}\put(76,20){\line(2,1){25}}
\put(25,10){\line(2,1){20}}\put(80,10){\line(2,1){20}}
\put(35,5){\line(2,1){10}}\put(90,5){\line(2,1){10}}
\put(145,80){\line(2,1){10}}\put(195,80){\line(2,1){10}}
\put(135,65){\line(2,1){25}}\put(196,70){\line(2,1){25}}
\put(145,60){\line(2,1){20}}\put(200,60){\line(2,1){20}}
\put(155,55){\line(2,1){10}}\put(210,55){\line(2,1){10}}
\put(280,70){\line(2,0){40}}
\put(265,80){\line(2,1){10}}\put(315,80){\line(2,1){10}}\put(375,80){\line(2,1){10}}
\put(255,65){\line(2,1){25}}\put(316,70){\line(2,1){25}}\put(376,70){\line(2,1){25}}
\put(265,60){\line(2,1){20}}\put(320,60){\line(2,1){20}}\put(380,60){\line(2,1){20}}
\put(275,55){\line(2,1){10}}\put(330,55){\line(2,1){10}}\put(390,55){\line(2,1){10}}
\put(60,25){\makebox(5,4){$e_1$}}
\put(25,22){\makebox(5,4){$\sigma_1$}}
\put(85,22){\makebox(5,4){$\sigma_2$}}
\put(60,85){\makebox(5,4){$e_2$}}
\put(25,82){\makebox(5,4){$\sigma_3$}}
\put(85,82){\makebox(5,4){$\sigma_4$}}
%\put(60,-10){\makebox(5,4){$Fig~~1$}}
 \put(150,30){\makebox(5,4){$e_1$}}
\put(210,30){\makebox(5,4){$e_2$}}
\put(145,72){\makebox(5,4){$\sigma_1$}}
\put(205,72){\makebox(5,4){$\sigma_2$}}
\put(180,0){\makebox(5,4){$Fig.~1$}}
\put(265,72){\makebox(5,4){$\sigma_1$}}
\put(325,72){\makebox(5,4){$\sigma_2$}}
\put(385,72){\makebox(5,4){$\sigma_3$}}
 \put(390,30){\makebox(5,4){$e_2$}}
  \put(300,75){\makebox(5,4){$e_1$}}
%\put(330,0){\makebox(5,4){$Fig~~3$}}
\end{picture}
\end{figure}
%\newpage
\indent Without loss of generality, we may suppose that $e_1\cap
e_2=\phi$,~and consider case (1).\\
\indent Let $e_1\in~E[\sigma_1,\sigma_2],~e_2\in
E[\sigma_3,\sigma_4]$,~and $C_i$~be the fundamental cycle
 in $T+e_i~(~1\leq i\leq 2~)$.

\vskip 0.2cm

\noindent {\bf Subcase A. }$C_1\cap C_2$ is a path.

\vskip 0.2cm \indent Let $P=C_1\cap C_2$ be a path with an end
vertex $x$ in $C_1\cap C_2$.~Let $e'_1$~and~$e'_2$~be two edges such
that $e'_1,~e'_2\in~E( T )$~,and $x\in~e'_1\cap~e'_2$.~Now
$e'_1,~e'_2\not\in~E(P)$.(~As shown in left hand side of
Fig.2).Consider a new spanning tree
$T'=T+\{e_1,~e_2\}-\{e'_1,~e'_2\}.$\\
\begin{figure}[h]
\centering
\begin{picture}(200,100)(0,0)
\thicklines \put(0,40){\circle{1.5}}\put(0,80){\circle{1.5}}
\put(50,30){\circle{1.5}}\put(50,80){\circle{1.5}}
\put(100,30){\circle{1.5}}\put(100,60){\circle{1.5}}
\put(30,100){\circle{1.5}}\put(70,100){\circle{1.5}}
 \put(0,40){\line(0,2){40}}\put(50,30){\line(0,5){50}}\put(100,30){\line(0,3){30}}
\put(30,100){\line(1,-1){20}}\put(50,80){\line(1,1){20}}
 \thinlines
\qbezier(0,40)(25,-10)(50,30)\qbezier(50,30)(80,-10)(100,30)
\qbezier(0,80)(10,110)(30,100)\qbezier(70,100)(100,110)(100,60)
\put(7,60){\makebox(5,4){$e_2$}}\put(90,45){\makebox(5,4){$e_1$}}
\put(55,50){\makebox(5,4){$P$}}\put(70,25){\makebox(5,4){$C_1$}}
\put(20,25){\makebox(5,4){$C_2$}}\put(45,73){\makebox(5,4){$x$}}
\put(65,85){\makebox(5,4){$e'_1$}}\put(40,95){\makebox(5,4){$e'_2$}}
\put(110,-5){\makebox(5,4){$Fig.2$}}
 \thicklines
\put(170,30){\circle{1.5}}\put(220,30){\circle{1.5}}\put(150,100){\circle{1.5}}
\put(220,60){\circle{1.5}}\put(170,80){\circle{1.5}}\put(170,60){\circle{1.5}}
\put(170,30){\line(0,5){50}}\put(220,30){\line(0,3){30}}\put(150,100){\line(1,-1){20}}
\put(170,80){\line(1,1){20}}
 \thinlines
\qbezier(120,40)(145,-10)(170,30)\qbezier(120,40)(120,100)(150,100)
\qbezier(170,30)(200,-10)(220,30)\qbezier(190,100)(220,110)(220,60)
\put(210,45){\makebox(5,4){$e_1$}}\put(160,80){\makebox(5,4){$x$}}
\put(175,50){\makebox(5,4){$P$}}\put(160,95){\makebox(5,4){$e_2$}}
\put(180,100){\makebox(5,4){$e'_1$}}\put(175,70){\makebox(5,4){$e'_2$}}
\put(190,25){\makebox(5,4){$C_1$}}\put(140,25){\makebox(5,4){$C_2$}}
%\put(170,-5){\makebox(5,4){$Fig~~5$}}
\end{picture}
\end{figure}

\vskip 0.2cm

 \noindent {\bf Subcase B.} $x\in~e_1$~or
$x\in~e_2$,~say $x\in~e_2$. ~(~As shown in right hand side of
Fig.2).

\vskip 0.1cm

If~$|E(P)| \geq 1$~,then we take edges $e'_1~\in C_1\backslash
E(T),~x\in e'_1,~e'_2\in E(P),~x\in e'_2~$. We construct a new
spanning tree~$T'=T+\{e_1,~e_2\}-\{e'_1,~e'_2\}$. If~$|E(P)|=0$,
then this may be a special case of A.

\vskip 0.2cm

\indent Let $T^{'}$ be the spanning tree as defined in either
subcase A or B. It is easy to see that $E(G_1)\backslash E(T')$~has
at most $\xi(G_0)$~odd components .~It
is contradictory to our suppose.~Therefore~$\xi(G_1)\leq\xi(G_0)$.\\
\indent Similarly, We may prove the claims in the cases of (2) and
(3).\\
\indent Repeat this procedure for~$G_2=G_1+\{e_3, e_4 \}, \cdots,
G_g=G_{g-1}+\{e_{2g-1}, e_{2g}\}$ until we get
$\xi(G_g)\leq\xi(G_{g-1})\leq\cdots\leq\xi(G_0)$, so
 $\gamma_{M}(G_g)\geq\gamma_{M}(G_0)+g=g$.\qed
\vskip 0.2cm
%\end{proof}
%\end{theorem}
%\begin{theorem}
\noindent{\bf Theorem 2.}{\it Let $G$ be a connected graph embedded
in an orientable surface $S_g$~and $T$~be a spanning tree of
$G$.~Then there are at least $2g$ noncontractible foundamental
cycles $C_1,~C_2,\cdots,C_{2g}$,~such that ~$C_{2i-1}\cap
C_{2i}\neq\phi $~for~
 $1\leq i\leq g$.~In particular, if $G$ is a
 one-face-embedded graph in $S_g$,~then for any spanning tree
 $T$~of $G$,~there are $2g$~edges in $G\backslash E(T)$~such
 that the corresponding $2g$ fundamental cycles $C_1,~C_2,\cdots,C_{2g}$~
 satisfy $C_{2i-1}\cap C_{2i}
\neq \phi$~for $1\leq i\leq g$.}
%\begin{proof}
\vskip 0.2cm

\noindent{\bf Proof}\quad. We contract $T$~into a single vertex
$v_T$~and delete all the possible edges on distinct faces. Then we
get a vertex-graph $G_T$~with exactly one vertex $v_T$~and one face
in $S_g$. There are two {\it crossed loops}, say
$e_\alpha,~e_\beta$,~such that the local rotation of semi-edges
incident to $v_T$~is $e_\alpha\cdots e_\beta\cdots e_\alpha\cdots
e_\beta$~( as shown in Fig.3). Furthermore, $e_\beta$ is the only
possible edge crossing $e_\alpha$ (~since otherwise $G_T$ would have
at least two faces! ) . Hence , all ( loop ) edges of $G_T$ may be
listed as follows: $e_1, e_2, \cdots, e_{2g-1}, e_{2g} $ such that
$e_{2i-1}$ crossing
  $e_{2i}$ for $i = 1, 2, \cdots, g$. It is easy to see that
  $e_{2i-1}$ and $e_{2i}$ determine two
fundamental cycles $C_{2i-1}$~and~$C_{2i}$~with a vertex in
common.\qed

\vskip 0.2cm

\begin{figure}[h]
\centering
\begin{picture}(80,100)(0,0)
\put(40,60){\circle*{2.5}}
 \put(10,90){\makebox(5,4){$e_\beta$}}
\put(65,90){\makebox(5,4){$e_\alpha$}}
\put(0,30){\makebox(5,4){$e_\alpha$}}
\put(70,30){\makebox(5,4){$e_\beta$}}\put(50,56){\makebox(5,4){$v_T$}}
\thicklines
 \qbezier(-10,80)(20,85)(40,60)
\qbezier(-10,40)(-40,60)(-10,80)
\qbezier(-10,40)(20,35)(40,60)\qbezier(90,40)(120,60)(90,80)
\qbezier(90,40)(60,35)(40,60) \qbezier(20,100)(15,85)(40,60)
\qbezier(60,100)(65,85)(40,60) \qbezier(20,100)(40,130)(60,100)
 \qbezier(20,20)(10,20)(40,60)\qbezier(90,80)(60,85)(40,60)
\qbezier(60,20)(65,20)(40,60) \qbezier(20,20)(40,0)(60,20)
\put(0,20){\line(1,1){80}}\put(0,100){\line(1,-1){80}}
 \thinlines
 \put(-10,50){\line(1,1){20}}\put(30,100){\line(1,1){10}}
 \put(28,90){\line(1,1){15}} \put(28,80){\line(1,1){20}}
  \put(40,70){\line(1,1){10}} \put(35,75){\line(1,1){10}}
\put(-20,60){\line(1,1){10}} \put(-15,55){\line(1,1){20}}
\put(0,45){\line(1,1){25}} \put(15,50){\line(1,1){10}}
\put(60,60){\line(1,1){10}} \put(70,50){\line(1,1){20}}
\put(80,45){\line(1,1){20}} \put(60,50){\line(1,1){20}}
\put(30,100){\line(1,1){10}}\put(30,20){\line(1,1){20}}
 \put(28,90){\line(1,1){15}} \put(28,80){\line(1,1){20}}
  \put(40,70){\line(1,1){10}} \put(35,75){\line(1,1){10}}
 \put(30,30){\line(1,1){15}}\put(32,40){\line(1,1){10}}
 \put(35,15){\line(1,1){15}} \put(45,15){\line(1,1){10}}
 \put(40,-5){\makebox(5,4){$Fig.3$}}
 \end{picture}
 \end{figure}
%\end{proof}
%\end{theorem}
\noindent{\bf Remark:}\quad Theorems 1 and 2 give a good
characterization of maximum genus of a graph(i.e., they implies the
existence of a polynomially bounded algorithm to find the maximum
genus of a graph).

\vskip 0.2cm

 \indent Let $T$ be a spanning tree
in $G$ with a group of fundamental cycles $C_1,C_2,\cdots,C_{2g}$.
If $C_{2i-1}\cap C_{2i} \neq \phi$ for $1\leq i\leq g$, then we say
 $\langle C_{2i-1},~C_{2i} \rangle$ is {\it an adjacent fundamental cycle
 pairs}
 $(1\leq i\leq g )$. If $g$ is chosen as the largest number
satisfying above condition, then we call $g$ the {\it maximum number
of adjacent fundamental cycle pairs} of $T$. Hence Theorem $2$
implies the following:

\vskip 0.2cm
%\begin{theorem}
\noindent{\bf Theorem 3}\quad{\it Any two spanning trees $T_1$ and
$T_2$ in a graph $G$ have the same maximum number of adjacent
fundamental cycle pairs. ( In fact , this unique number is
$\gamma_{M}(G)$, the maximum genus of $G$ )}.
%\end{theorem}
\vskip 0.2cm

 This generalizes a result of Fu et al[2] where they
introduced the concept {\it intersecting graph} which is determined
by bases of cycle space of a graph to describe the maximum genus of
a graph. In fact, our result stands for any spanning tree's
fundamental cycles.

\vskip 0.2cm

\newtheorem{corollary}{Corollary}
\begin{corollary}
~If a connected graph $G$~has a spanning tree $T$~such that any two
fundamental cycles have a vertex in common.~Then $G$~is
upper-embeddable.
\end{corollary}

\vskip 0.2cm

 Sometimes however, we need a refined form of Theorems 1
and 2 in practice. The following result gives us a recursive
relation between the maximum genera of a graph and its subgraph(s).

\vskip 0.2cm
%\begin{theorem}
\noindent{\bf Theorem 4}\quad{\it Let $G$ be a connected graph and
$T$ be an arbitrary spanning tree in $G$. If $e_1, e_2$ are two
edges not in $G$ and the two cycles $C_T(e_1)$ and $C_T(e_2)$ have a
vertex in common. Then $\gamma_M(G)=\gamma_M(G+e_1+e_2)-1$. In
particular, $G$ is upper-embeddable  if and only if $G+e_1+e_2$ is
upper-embeddable.}
%\end{theorem}
\vskip 0.2cm

One may easily see that this generalizes a recursive relation for
maximum genus of Xuong[9] and ( we will see in the next section )is
much more practical in use.

\section{Applications }
\indent Now in this section, we begin to apply Theorems $1-2$ to
determine the maximum genus of some type of graphs.

Let us recall that the essence of Xuong's method[9] consists of two
parts: one is to find an optimal tree $T$ in a graph $G$ having the
smallest number of odd components; the other is to organize edges of
$E(G)\backslash E(T)$ into adjacent pairs such as
$$
E(G)\backslash E(T)=\{e_1, e_2,
\cdots,~e_{2s}\}\cup\{~f_1,~f_2,~\cdots~f_m\},
$$
where $e_{2i-1}\cap~e_{2i}\neq\phi$($1\leq~i\leq~s$~)~and~
$C_T(f_i)\cap{C_T(f_j)}=\phi$,~for~$1\leq~i<j\leq~m$~and
~$s=\gamma_{M}(G),~m=\xi(G)$. Compared with the above procedure,
Theorems 1 and 2 consider adjacent foundamental cycle pairs(rather
than adjacent pairs of edges). We may construct large genus
embedding from any spanning tree $T$, although it may have very
large number of odd components in $G\backslash E(T)$. This greatly
releases the conditions of Xuong. Of course, an optimal tree is also
valid in our constructions. Hence, Theorems 1 and 2 generalize
Xuong's characterization of maximum genus. Based on this idea, we
may construct a large orientable genus as follows: Take a specific
spanning tree $T$ in graph $G$ and first organize some non-tree
edges into adjacent pairs ( as Xuong did ) and then match other
possible non-tree edges into pairs such that their fundamental
cycles also become adjacent fundamental cycle pairs. It is easy for
one to see that the second part of non-tree edges may be chosen as
an edge-cut of $G$. Therefor, Theorems $1-2$ may be useful in
determination of a maximum genus of a graph $G$ with a specific
edge-cut. Now, the following result is easy to be verified.

\vskip 0.2cm

%\begin{theorem}
\noindent{\bf Theorem 5}\quad{\it Let $A=\{e_1,~e_2,\cdots,~e_k~\}$
be an edge-cut of $G$ such that $G-A$ has exactly two components
$G_1$ and $G_2$. If both $G_1$ and $G_2$ are upper-embeddable, then
$~\gamma_{M}(G)\geq\lfloor\frac{\beta(G)}{2}\rfloor-1$ .
Furthermore,if $G$ satisfies one of the
following conditions , then $G$ is upper-embeddable:\\
(1). $\beta(G_1)\equiv\beta(G_2)\equiv~0(~mod~2~)$\\
(2). $|A|\equiv~1(~mod~2~)$ and
$\beta(G_1)+\beta(G_2)\equiv~1(~mod~2~)$.}
%\end{theorem}
\vskip 0.2cm

 \indent The next result is due to Huang. As a
consequence of the above results, we will give another proof.

\vskip 0.2cm

%\noindent{\bf Theorem 6(Huang[5])}

%\begin{theorem}{\upshape{(Huang[5])}}
\noindent{\bf Theorem 6(Huang[5])}\quad{\it Let $G$~be a strongly
embedded graph in an orientable surface $S_g$~( i.e., all facial
walks are cycles ).~If the dual graph $G^*$~of $G$~has a surface
separating Hamiltonian cycle, then $G$~is upper-embeddable.}
%\begin{proof}
\vskip 0.2cm

 \noindent{\bf Proof}\, We will show the existence of a
spanning tree $T$~of $G$~satisfying the conditions in Theorems 1 and
2.~Let $\mathcal{F}=\{f_1,~f_2,\cdots~f_\varphi\}$~be the face-set
of $G$~and $C^*$~be a surface separating Hamiltonian cycle in
$G^*$.~Let $E(C^*)=\{e^*_1,~e^*_2,\cdots~e^*_\varphi~\}$~and
$e^*_i~=(f_i,f_{i+1})$~for $1\leq i\leq\varphi$.~Let $e_i$~be the
edge in $\partial f_i\cap\partial f_{i+1}$~corresponding to
$e^*_i$~for $1 \leq i\leq\varphi$ (where
$\partial f_i$~denotes the boundary of $f_i$).\\
{\bf Claim 2.}~$G-\{e_1,~e_2,\cdots~e_{\varphi-1}\}$~is a one-face
embedded subgraph of $G$~in $S_g$,~Furthermore, $G-
\{e_1,~e_2,\cdots~e_\varphi\}$~has exactly two components
$G_1$~and $G_2$.\\
\indent  Now $E[G_1,~G_2]=\{e_1,~e_2,\cdots~e_\varphi\}$. Let
$G_1\subset~Int(C^*)$~and $G_2\subset~Ext(C^*)$~and $\partial
f_i$~denotes  the boundary cycle of $f_i$ for $1\leq i\leq
\varphi$.~Then we may construct a graph as follows.~ $H_0=(\partial
f_1\cup\partial f_2\cup\cdots\cup\partial f_{\varphi
-1})\setminus\{e_1,~e_2,\cdots~e_{\varphi-1}\}$.~It is easy to see
that $H_0$ is a connected spanning subgraph of $G-
\{e_1,~e_2,\cdots~e_{\varphi-1}\}$.~(~Hence,~a spanning subgraph of
$G$).~Let $e_\varphi=(\alpha, \beta)$~with
$\alpha\in~V(G_1),~\beta\in~V(G_2)$.~Then $H_0-e_\varphi$~has
exactly two components $H',~H_1$~with
$H'=G_1$.\\
\vskip 0.2cm

 {\bf Claim 3.}~If $H_1$~has a cycle $C$,~then $C$~must
be a noncontractible cycle.\\
\indent This follows from the fact that $S_g-H_0$~has only one
component. If $H_1$~has a cycle $C_1$,~then delete an edge $e'_1
\in~C_1$~and get a subgraph $H_2$~of $H_1$~with $V(H_2)=V(H_1
)$.~Repeat this procedure until we arrive at a connected subgraph
$H_k$~of $H_1$~with $V(H_k)=V(H_1)$~and $H_k$~has no cycle.\\

{\bf Claim 4:}~$T=H'\cup~H_k\cup~\{e_\varphi\}$~is a spanning tree
of $G$,~such that each fundamental cycle $C_T(e_i)$~in $T+e_i$~has
an edge $e_\varphi$~in common for $i=1, 2,\cdots,\varphi-1~$.\\
\indent To see this, we consider an edge $e_i=(x_i,y_i)\in[ H'~,~H_k
]\subseteq[G_1,~G_2]$,~such that $x_i\in H',~y_i\in H_k$.~Since $H'(
H_k )$~is connected, there is a path$P_i(Q_i)$~in $H'(H_k)$~joining
$\alpha(\beta)$~and $x_i (y_i)$.~Hence, $C_T(e_i
)=\{e_\varphi\}\cup~P_i\cup~Q_i\cup~\{e_i\}$~is a cycle
containing $e_\varphi$~for $1\leq i\leq\varphi$.\\
Now we find a spanning tree $T$~of $G$~such that:\\
$(1)$. All the fundamental cycles
$C_T(e_1),~C_T(e_2),\cdots,~C_T(e_{\varphi-1})$
~has an edge in common;\\
$(2)$. By Theorems 1 and 2, and the fact that $T$~is also a spanning
tree in $G-\{e_1,~e_2,\cdots~e_{\varphi-1}\}$, there are another
group of fundamental ( noncontractible ) cycles
$C_1,~C_2,\cdots,~C_{2g}$~such that
$C_{2i-1}\cap~C_{2i}\neq{\o}$~for $1\leq i\leq g$.~ By theorem 1,
$G$~is upper-embeddable.\qed

%\end{proof}
%\end{theorem}

\vskip 0.2cm

 One may readily see that a surface separating cycle may
not be Hamiltonian and the hosting surface on which a graph is
embedded may not be orientable. Thus, Theorem 6 can be extended to a
much more generalized form.

\vskip 0.2cm
%\begin{theorem}
\noindent{\bf Theorem 7}\quad{\it Let $G$ be an embedded graph in a
surface $\sum$ such that the dual graph $G^*$ of $G$ contains a
surface separating cycle $C^*$ such that both of the left subgraph
$G_L(C^*)$ and right subgraph $G_R(C^*)$ of $G$ are upper embeddable
. Then $~\gamma_{M}(G)\geq\lfloor\frac{\beta(G)}{2}\rfloor-1$. In
particular, if ${\beta(G_L(C^*))}\equiv {\beta(G_R(C^*))}\equiv
0\,(\,mod\, 2)$, then $G$ is upper-embeddable.}

%\end{theorem}

\vskip 0.2cm

 \noindent{\bf Remark:}\, The term ``~left( right )
subgraph~'' follows from [7]
\begin{corollary}
~~If $G$ is an embedded graph on the Klein bottle such that the dual
graph $G^*$ has a surface separating Hamiltonian cycle. Hence
$\gamma_{M}(G)\geq\lfloor\frac{\beta(G)}{2}\rfloor-1$.
\end{corollary}
\vskip 0.2cm

 \indent In practical use, our attentions need not to be
restricted to graphs with an edge-cut. Theorems 1-4 provide us a
tool to evaluating large genus embeddings in more extended range of
graphs. The following results show us how to do so ( we omit the
proof of them ).
\\
\vskip 0.1cm

\par
%\begin{theorem}
\noindent{\bf Theorem 8.}\quad{\it
The following graphs are upper-embeddable:\\
\noindent(1). The cartisian product $G\times P_n$ of a simple
connected graph $G$ and a path $P_n$ with $n(\geq 1)$ egdes;\\
 \noindent (2). The
composition of two disjoint Halin graphs $H_1$ and $H_2$ with some
edges $e_1, e_2,\cdots,e_k (k\geq 2)$ connecting
them;\\
\noindent (3). The $n$-cube $Q_n$ which is composed of two
$(n-1)$-cube $Q_{n -1}$ together with some edges joining the two
copies of vertice in $Q_{n-1}$.\\
\noindent (4). The generalized Petersen graphs $P(n , k)$ which is
determined by $n-$cycle $(u_1,u_2,\cdots,u_n)$ and vertices
$v_1,v_2,\cdots,v_n$ such that (i) each $(u_i,v_i)\in E, 1\leq i\leq
n$;(ii) $(u_i, v_{i+k})\in E, 1\leq i\leq n$.}

\vskip 0.2cm
%\begin{proof}
%\end{theorem}
Note: A graph $G=(V,~E)$~is a {\it Halin graph} if~$G$~is obtained
by joining the leaves(1-valent vertices) of a plane tree $T$ with a
cycle in this orientation and the definition of cartisian product of
two graphs may be fund in and textbook of graph theory.

\vskip 0.3cm

%\newpage
\section{ A polynomially bounded algorithm}

In this section we shall present a polynomially bounded algorithm to
find the maximum genus of a given graph. A basic fact is that
Theorems 1 and 2 present a {\it good characterization } of maximum
genus problem, i.e., we have the following

\vskip 0.2cm
 \noindent{\bf Theorem 9}\quad {\it To determine the
maximum genus of a graph $G$ is equivalent to determine a maximum
matching of the graph $G_M=(V_M,E_M)$, called {\it fundamental
intersecting graph of $G$}, where $V_M$ is the set of fundamental
cycles of a spanning tree $T$ of $G$ and any two cycles in $V_M$ are
adjacent if and only if they have at least a vertex in common.}

\vskip 0.2cm

 We observe that the fastest algorithm to find a maximum
matching in a graph $G$ is due to Micali-Vazirani[8] which will end
in $O(m\sqrt{n})$ steps, where $m$ and $n$ are, respectively, the
number of edges and vertices of $G$. Based on this fact and Theorems
1 and 2 we may construct a new algorithm to determine the maximum
genus of a graph $G$.

\vskip 0.2cm

\noindent{\bf Fundamental cycle algorithm}

\vskip 0.2cm \noindent{\bf Step 1.Input the date of the graph $G$
and then searching for a spanning tree $T$ and the set $V_M$ of fundamental cycles in $G$;}\\
 \noindent{\bf Step 2. For cycles in $V_M$ we build the graph
 $G_M$};\\
 \noindent{\bf Step 3 Perform Micali-Vazirani algorithm to find a maximum matching in
 $G_M$ and then terminate.}

 \noindent{\bf Remark:} Since the number of fundamental cycles in a graph $G$ of order $n$
 is ${(\beta(G))}$, this algorithm will end in at most $O({(\beta(G))}^{\frac{5}{2}})$
 steps. Although Furst, Gross and McGeoch had already construct the
 first polynomially bounded algorithm[3], this result is a new approach
 to do so.

\end{document}